\documentclass[11pt,leqno]{article} 
\usepackage{graphics}
\newtheorem{thm}{Theorem}[section]
\newtheorem{lma}{Lemma}[section]
\newtheorem{cor}{Corollary}

\newcommand{\beqa}{\begin{eqnarray}}
\newcommand{\eeqa}{\end{eqnarray}}

\newcommand{\pf}{\noindent {\bf Proof:} $\s$ }
\newcommand{\epf}{ \hfill$\diamondsuit$ \medskip}
\newcommand{\md}{\medskip}

\newcommand{\beq}{\begin{equation}}
\newcommand{\eeq}{\end{equation}}
\newcommand{\lbl}{\label}
\newcommand{\s}{\; \;}

\newcommand{\ep}{\epsilon}

\newcommand{\la}{\lambda}

\newcommand{\ra}{\rightarrow}

\newcommand{\p}{\varphi}

\title{Global solution curves for several  classes of singular periodic problems}

\author{
Philip Korman   \\ 
Department of Mathematical Sciences \\ 
University of Cincinnati \\ 
Cincinnati Ohio 45221-0025 \\
}

\date{}

\begin{document}

\maketitle
\begin{abstract} 
Using  continuation methods and bifurcation theory, we study the exact multiplicity of periodic solutions, and the global solution structure, for three  classes of periodically forced  equations with singularities, including the equations arising in micro-electro-mechanical systems (MEMS), the ones in  condensed matter physics, as well as  A.C. Lazer and S. Solimini's \cite{LS} problem. 
 \end{abstract}

\begin{flushleft}
Key words: Periodic solutions, the  exact multiplicity of solutions. 
\end{flushleft}

\begin{flushleft}
AMS subject classification:  34C25, 34B16, 78A35, 74F10.
\end{flushleft}

\section{Introduction}
\setcounter{equation}{0}
\setcounter{thm}{0}
\setcounter{lma}{0}

The study of periodic solutions of equations with singularities began with a remarkable paper of A.C. Lazer and S. Solimini \cite{LS}. One of the model problems considered in that paper involved (in case  $c=0$)
\beq
\lbl{i1}
u''(t)+c u'(t)+\frac{1}{u^p(t)}=\mu +e(t), \s\s u(t+T)=u(t) \,,
\eeq
where $e(t)$ is a continuous function of period $T$ with  $\int_0^T e(t) \, dt=0$, $c \in R$, and $p>0$. (Here and throughout the paper we write $u(t+T)=u(t)$ to indicate that we are searching for  $T$-periodic solutions.) A.C. Lazer and S. Solimini \cite{LS} proved that the problem (\ref{i1}) has a positive $T$-periodic solution if and only if $\mu >0$. It turned out later that problems with singularities occur often in applications. The recent book of P.J. Torres \cite{T} contains a review of these applications, with up to date results, and a long list of open problems.
\medskip

We approach the problem (\ref{i1}) with continuation (bifurcation) methods, and add some details to the Lazer-Solimini result, which turn out to be useful, in particular, for numerical computation of solutions. Write $u(t)=\xi+U(t)$, with $\int_0^T U(t) \, dt=0$, so that $\xi$ is the average of the solution $u(t)$. We show that $\xi$ is a {\em global parameter} i.e., the value of $\xi$ uniquely identifies both $\mu$ and the corresponding $T$-periodic solution $u(t)$ of (\ref{i1}). Hence, the set of positive $T$-periodic solutions of (\ref{i1}) can be represented by a curve in $(\xi, \mu)$ plane. We show that this curve, $\mu =\p (\xi)$, is hyperbola-like, where  $\p (\xi)$ is a decreasing function, defined for all $\xi \in (0,\infty)$, and $\lim _{\xi \ra 0} \p (\xi)=\infty$, $\lim _{\xi \ra \infty} \p (\xi)=0$. 
It turns out to be relatively easy to compute numerically the solution curve $\mu =\p (\xi)$, see Figure 1 below. In the last section we  explain the implementation  of our  numerical  computations, using the {\em Mathematica} software. Using recent results of R. Hakl and M. Zamora \cite{H}, we also discuss the solution curve in case $e(t) \in L^2$.
\medskip

There is  a considerable recent interest in micro-electro-mechanical systems (MEMS), see e.g., J.A. Pelesko \cite{P}, Z. Guo and J. Wei \cite{GW}, or P. Korman \cite{K1}.  Recently A. Guti\'{e}rrez and P.J. Torres \cite{G} have considered an idealized mass-spring model for MEMS, which they reduced to the following problem
\beq
\lbl{i2}
u''(t)+c u'(t)+bu(t)+\frac{a(t)}{u^p(t)}=\mu +e(t), \s\s u(t+T)=u(t) \,,
\eeq
with $p=2$, and $e(t)=0$ ($b>0$ and $c \in R$ are constants). We show that the solution curve $\mu =\p (\xi)$ of (\ref{i2}) is parabola-like, and there is a $\mu _0>0$ so that the problem (\ref{i2}) has exactly two positive $T$-periodic solutions for $\mu>\mu _0$, exactly one positive $T$-periodic solution for $\mu=\mu _0$, and no positive $T$-periodic solutions for $\mu<\mu _0$, see Figure $2$. This extends the corresponding result in \cite{G}.
\medskip

For both of the above equations  we were aided by the fact that the nonlinearities were convex.
An interesting problem where the nonlinearity changes convexity arises as a model for fluid adsorption and wetting on a periodic corrugated substrate:
\beq
\lbl{i3}
\s u''(t)+c u'(t)+a\left(\frac{1}{u^4(t)}-\frac{1}{u^3(t)} \right)=\mu +e(t), \s\s u(t+T)=u(t) \,,
\eeq
in case $c=0$,  and $\mu =0$, see C. Rasc\'{o}n et al \cite{R}, and also  P.J. Torres \cite{T}, where this problem was suggested as an open problem. Using an idea from G. Tarantello \cite{Ta}, we again describe the shape of the solution curve, and obtain an exact multiplicity result, see Figure $3$. For the physically significant case when $c=0$,  and $\mu =0$, our result implies the existence  and uniqueness of positive $T$-periodic solution.
\medskip

Remarkably, in all three cases the graph of  the global solution curve is similar to that of the nonlinearity $g(u)$.
\medskip

We now outline our approach. 
We embed (\ref{i1}) into a family of problems
\beq
\lbl{i4}
u''(t)+c u'(t)+k\frac{1}{u^p(t)}=\mu +e(t), \s\s u(t+T)=u(t) \,,
\eeq
with the parameter $0 \leq k \leq 1$. 
When $k=0$ and $\mu =0$, the problem is linear, and it has a unique $T$-periodic solution of any average  $\xi$. We show that if the average of solutions is kept fixed, the solutions $(u,\mu)$ of (\ref{i4}) can be continued, using the implicit function theorem,  for all $0 \leq k \leq 1$, and at $k=1$ we obtain a unique solution  $(u,\mu)$ of (\ref{i2}) of any average  $\xi$. It follows that the value of $\xi$ is a global parameter. We then obtain the exact shape of  the solution curve $\mu=\phi(\xi )$, by calculating the sign of $\mu ''(\xi _0)$ at any critical point $\xi _0$, or by proving that $\phi(\xi )$ is monotone.
\medskip

\section{Preliminary results}
\setcounter{equation}{0}
\setcounter{thm}{0}
\setcounter{lma}{0}

We consider $T$-periodic functions, and use $\omega=\frac{2 \pi}{T}$ to denote the frequency. We denote by $L^2_T$ the subspace of $L^2(R)$, consisting of $T$-periodic functions. We denote by $H^2_T$ the subspace of the Sobolev space $H^2(R)=W^{2,2}(R)$, consisting of $T$-periodic functions. By $\bar L^2_T$ and $\bar H^2_T$ we denote the respective subspaces of $L^2_T$ and $H^2_T$, consisting of  functions of zero average on $(0,T)$, i.e., $\int_0^T u(t) \, dt=0$. 
The following lemma we proved in \cite{K}, by using Fourier series.
\begin{lma}\lbl{lma:1}
Consider the linear problem 
\beq
\lbl{1}
y''(t)+\la y'(t)=e(t),
\eeq
with $e(t) \in \bar L^2_T$ a given function of period $T$, of zero average, i.e.,  $\int_0^T e(s) \, ds=0$. Then the problem (\ref{1}) has a unique $T$-periodic solution $u(t) \in H^2_T$ of any average.
\end{lma}

The following lemma is well-known as Wirtinger's inequality. Its proof follows easily by using the complex Fourier series, and the orthogonality of the functions $\{e^{i \omega n t} \}$ on the interval $(0,T)$.
\begin{lma}\lbl{lma:2}
Assume that $f(t)$ is a continuously differentiable function  of period $T$, and of zero average, i.e.,  $\int_0^T f(s) \, ds=0$. Then
\[
\int _0^T {f'}^2(t) \, dt \geq \omega ^2 \int _0^T f^2(t) \, dt.
\]
\end{lma}

We consider next the following linear periodic problem in the class of functions of zero average 
\beq
\lbl{2} \hspace{0.2in}
w''(t)+c w'(t)+h(t)w(t)=\mu, \s w(t+T)=w(t), \s \int_0^T w(s) \, ds=0,
\eeq
where $h(t)$ is a given continuous function of period $T$,  $\mu$ is a parameter.
\begin{lma}\lbl{lma:3}
Assume that 
\beq
\lbl{3}
h(t) <  \omega^2.
\eeq
Then the only solution of (\ref{2}) is $\mu=0$ and $w(t) \equiv 0$.
\end{lma}

\pf
Multiplying the equation (\ref{2}) by $w$, and integrating
\[
\omega ^2 \int_0^T w^2 \, dt \leq \int_0^T {w'}^2 \, dt=\int_0^T h(t)w^2 \, dt < \omega ^2 \int_0^T w^2 \, dt \,.
\]
It follows that $w(t) \equiv 0$, and then $\mu =0$.
\epf

\noindent
{\bf Remark} A similar result holds if (\ref{3}) is replaced by $|h(t)| \leq \omega \sqrt{\omega ^2+c^2}$, see J. Cepicka et al \cite{D}, or Lemma 2.3 in \cite{K}. However, for singular problems we need a one-sided condition (\ref{3}).
\medskip

We consider next  another linear problem (now $w(t)$ is not assumed to be of zero average)
\beq
\lbl{5} 
w''(t)+c w'(t)+h(t)w(t)=0, \s\s w(t+T)=w(t).
\eeq
Here again $h(t)$ is a given continuous function of period $T$.
\begin{lma}\lbl{lma:4}
Assume that for all $t$
\beq
\lbl{6}
h(t)< \frac{c^2}{4}+\omega^2 \,.
\eeq
Then any non-trivial solution of (\ref{5}) is of one sign, i.e., we may assume that $w(t)>0$ for all $t$.
\end{lma}

\pf
Assume, on the contrary, that we have a sign changing solution $w(t)$. Since $w(t)$ is $T$-periodic, there exist $t_1<t_2$, such that $t_2-t_1=T$, and
$
w(t_1)=w(t_2)=0$.
Since we also have $w'(t_1)=w'(t_2)$, it follows that $w(t)$ has at least one more root on $(t_1,t_2)$. If we now consider   an eigenvalue problem on $(t_1,t_2)$ 
\beq
\lbl{7} \hspace{0.3in}
w''(t)+c w'(t)+\la h(t)w(t)=0, \s \mbox{for $t_1 <t <t_2$}, \s\s w(t_1)=w(t_2)=0,
\eeq
it follows that $w(t)$ is the second or higher eigenfunction, with the corresponding eigenvalue $\la =1$  being the second or higher eigenvalue, so that $\la _2 \leq 1$. On the other hand, we consider the following  eigenvalue problem on $(t_1,t_2)$
\beq
\lbl{8}
z''(t)+c z'(t)+\nu \left(\frac{c^2}{4}+\omega^2 \right)z(t)=0, \s\s z(t_1)=z(t_2)=0.
\eeq
Its eigenvalues  $\nu _n=\frac{\frac{c^2}{4}+\frac{n^2\omega^2}{4}}{ \frac{c^2}{4}+\omega^2}$ are smaller than the corresponding eigenvalues of (\ref{7}) (see e.g., p.$174$ in \cite{W}). We then have $\nu _2=1<\la _2 \leq 1$, a contradiction.
\epf

We shall also need the adjoint linear problem
\beq
\lbl{9} \hspace{0.3in}
v''(t)-c v'(t)+h(t)v(t)=0,  \s\s v(t+T)=v(t)=0.
\eeq
\begin{lma}\lbl{lma:5}
Assume the condition (\ref{6}) holds. If the problem (\ref{5}) has a non-trivial solution $w(t)$, the same is true for the adjoint problem (\ref{9}). Moreover, we then have $v (t)>0$ for all $t$.
\end{lma}

\pf
Assume that the problem (\ref{5}) has a non-trivial solution, but the problem (\ref{9}) does not. The differential operator given by the left hand side of (\ref{9}) is Fredholm, of index zero. Since its kernel is empty, the same is true for its co-kernel, i.e. we can find a solution $z(t)$ of 
\beq
\lbl{10} \hspace{0.3in}
z''(t)-c z'(t)+h(t)z(t)=w(t),  \s\s z(t+T)=z(t)=0.
\eeq
Multiplying the equation (\ref{10}) by $w(t)$, the equation (\ref{5}) by $z(t)$, subtracting and integrating, we have
\[
\int_0^T w^2(t) \, dt=0,
\]
a contradiction.
\md

Positivity of $v (t)$ follows by  the previous Lemma \ref{lma:4} (in which no assumption on the sign of $c$ was made).
\epf
\begin{lma}\lbl{lma:6}
Consider the non-linear problem
\beq 
\lbl{11d}
u''(t)+c u'(t)+g(u(t))=\mu +e(t), \s\s u(t+T)=u(t), 
\eeq
with $c,\mu \in R$, $e(t) \in \bar L^2_T[0,T]$, $g(u) \in C^1(R)$. 
Assume that for some $\omega _1 >0$
\beq
\lbl{11e}
g'(u) \leq \omega _1 <  \omega^2 \,, \s \mbox{for all $u \in R$} \,.
\eeq
Then there is a constant $c_0>0$, so that any solution of (\ref{11d}) satisfies
\beq
\lbl{11l}
\int_0^T \left({u''}^2(t)+{u'}^2(t) \right)\, dt \leq c_0 \,, \s \mbox{uniformly in $\mu$} \,.
\eeq
\end{lma}

\pf
Multiply the equation (\ref{11d}) by $u''$ and integrate
\[
\int_0^T {u''}^2(t) \, dt-\int_0^T g'(u){u'}^2(t) \, dt = \int_0^T u''(t) e(t) \, dt \,.
\]
Using the estimate $\int_0^T g'(u){u'}^2(t) \, dt \leq   \omega _1  \int_0^T {u'}^2(t) \, dt \leq \frac{\omega _1}{\omega^2} \int_0^T {u''}^2(t) \, dt$,
we get a bound on $\int_0^T {u''}^2(t) \, dt$, and using Wirtinger's inequality, we complete the proof.
\epf
\begin{cor}
Write any solution of (\ref{11d}) as $u=\xi+U$, with $\int_0^T U(t) \, dt=0$.
Then we have the estimate
\[
||U||_{H^2} \leq c_0 \,, \s \mbox{uniformly in $\xi$ and $\mu$} \,.
\]
\end{cor}

\pf
In (\ref{11l}), we have $u'=U'$, $u''=U''$, then use Wirtinger's inequality to estimate $\int_0^T {U}^2(t) \, dt$.
\epf

The following lemma is easily proved  by integration.
\begin{lma}\lbl{lma:i6}
Assume that any $T$-periodic solution  $u(t)$ of the problem (\ref{11d}) satisfies $|g(u(t))| \leq M, \s \mbox{for all $t \in R$}$. Then
\[
|\mu | \leq M.
\]
\end{lma}

\section{Continuation of solutions}
\setcounter{equation}{0}
\setcounter{thm}{0}
\setcounter{lma}{0}
We begin with continuation of solutions of any fixed average.
We consider the problem
\beq 
\lbl{11}
u''(t)+c u'(t)+kg(u(t))=\mu +e(t), \s\s u(t+T)=u(t), 
\eeq
\beq
\lbl{11a}
\frac{1}{T} \int_0^T u(t) \, dt=\xi.
\eeq
where $0 \leq k \leq 1$, $\mu$, $c$ and $\xi $ are parameters, and  $e(t) \in \bar L^2_T$.
\begin{thm}\lbl{thm:1}
 Assume that the function $g(u) \in C^1(R)$ satisfies the condition (\ref{11e}), and any solution of (\ref{11}-\ref{11a}) satisfies $|g(u(t))| \leq M, \s \mbox{for all $t \in R$}$, and all $0 \leq k \leq 1$, with some constant $M$.
Then  one can find a unique $\mu=\mu (k,\xi)$ for which the problem (\ref{11}), (\ref{11a}) has a unique classical solution. (I.e., for each $\xi$ there is a unique solution pair $(\mu, u(t))$.)
\end{thm}

\pf
We begin by assuming that $\xi=0$. We wish to prove that there is a unique $\mu=\mu (k)$ for which the problem (\ref{11}) has a  solution of zero average, and that solution is unique. We recast the equation (\ref{11}) in the operator form
\beq
\lbl{14}
F(u,\mu,k)=e(t), 
\eeq
where $F \; :   \bar H^2_T \times R \times R \ra L^2_T$ is defined by
\[
F(u,\mu,k)=u''(t)+c u'(t)+kg(u(t))-\mu.
\]
When $k=0$ and $\mu=0$, the problem (\ref{14}) has a unique $T$-periodic solution of zero average, according to the Lemma \ref{lma:1}. We now continue this solution for  increasing $k$, i.e., we solve (\ref{14}) for the pair $(u,\mu)$ as a function of $k$. Compute the Frechet derivative
\[
F_{(u,\mu)}(u,\mu,k)(w, \mu ^*)=w''(t)+c w'(t)+kg'(u(t))w(t)-\mu^*.
\]
According to Lemma \ref{3} the only solution of the linearized equation
\[
F_{(u,\mu)}(u,\mu,k)(w, \mu ^*)=0, \s\s w(t+T)=w(t)
\]
is $(w, \mu ^*)=(0,0)$. Hence, the map $F_{(u,\mu)}$ is injective. Since this map is Fredholm of index zero, it is also surjective. The implicit function theorem applies, giving us  locally  a curve of solutions $u=u(k)$ and $\mu=\mu(k)$. This curve continues for all $0 \leq k \leq 1$ because  the curve cannot go to infinity at some $k$, by Lemmas \ref{lma:6} and \ref{lma:i6}. (We have the uniform boundness of the solution by the Sobolev embedding theorem, and then we bootstrap to the boundness in $C^2$, since $\mu$ is bounded by Lemma \ref{lma:i6}. Therefore solutions of class $H^2_T$ are in fact classical.)
\md

The  solution pair $(u,\mu)$, which we found at the parameter value of $k=1$, is unique since otherwise we could continue another solution pair  for decreasing $k$, obtaining   at $k=0$ a  $T$-periodic solution of zero average for 
\[
u''+c u'=\mu _0 +e(t),
\]
with some constant $\mu_0$. By integration, $\mu_0 =0$, and then we obtain another solution of the problem (\ref{1})  of zero average  (since curves of solutions of zero average  do not intersect, in view of the implicit function theorem), contradicting Lemma \ref{lma:1}.
\md

Turning to the solutions of any average $\xi$, we again solve (\ref{14}), but redefine $F \;\; :  \bar H^2_T \times R \times R \ra L^2_T$ as follows
\[
F(u,\mu,k)=u''(t)+c u'(t)+kg(u(t)+\xi)-\mu.
\]
As before, we obtain a solution $(u,\mu)$  at $k$, with $u$ of zero average, which implies that $(u+\xi,\mu)$ is a solution of  our problem (\ref{11}), (\ref{11a}) of average $\xi$.
\epf

Next we discuss the continuation in $\xi$ for fixed $k$. As before, solutions of class $H^2_T$ are in fact classical.

\begin{thm}\lbl{thm:11}
Assume that the condition (\ref{11e}) holds for the problem (\ref{11}). Then solutions of (\ref{11}) can be locally continued in $\xi $. We have a continuous solution curve $(u,\mu)(\xi) \subset H^2_T \times R$, with $\xi$ being a global parameter (i.e., the value of $\xi$ uniquely identifies the solution pair $(u(t),\mu)$).
\end{thm}

\pf
The proof is essentially the same as for continuation in $k$ above. Writing $u(t)=\xi+v(t)$, with $v \in \bar H^2_T$, we recast the equation (\ref{11}) in the operator form
\[
F(v,\mu,\xi)=e(t), 
\]
where $F \; :   \bar H^2_T \times R \times R \ra L^2_T$ is defined by $F(v,\mu,k)=v''(t)+c v'(t)+kg(\xi+v)-\mu$.
We show exactly as before that the implicit function theorem applies, allowing local continuation of solutions.
\epf

\section{The global solution curve for Lazer-Solimini problem}
\setcounter{equation}{0}
\setcounter{thm}{0}
\setcounter{lma}{0}

We consider positive $T$-periodic solutions of  ($p>0$, $c \in R$)
\beq
\lbl{l1}
u''(t)+c u'(t)+\frac{1}{u^p(t)}=\mu +e(t), \s\s u(t+T)=u(t) \,.
\eeq

\begin{lma}\lbl{lma:l0}
Assume that $e(t) \in C(R)$ is T-periodic, with $\int_0^T e(t) \, dt=0$, $p>0$. There is an $\ep =\ep(\mu)>0$, so that any solution of (\ref{l1}) satisfies $u(t)>\ep$ for all $t$.
\end{lma}

\pf
Let $m=\min _t u(t)$, and $u(t_0)=m$. Evaluating the equation (\ref{l1}) at $t_0$, we have 
\[
\frac{1}{m^p} \leq \mu +e(t_0) \leq \mu +\max _te(t)\,,
\]
and the proof follows.
\epf

\begin{thm}\lbl{thm:l1}
Assume that the function $e(t) \in C(R)$ is T-periodic, with $\int_0^T e(t) \, dt=0$, $p>0$, and $c \in R$. The problem (\ref{l1}) has a  unique positive $T$-periodic classical solution if and only if $\mu >0$. Moreover, if we decompose $u(t)=\xi+U(t)$ with $\int_0^T U(t) \, dt=0$,   then the value of $\xi$ is a global parameter, uniquely identifying the solution pair $(u(t),\mu)$, and  all positive $T$-periodic solutions lie on a unique hyperbola-like curve $\mu=\p (\xi)$. The function $\p (\xi)$ is continuous and decreasing, it is defined for all $\xi \in (0,\infty)$, and $\lim _{\xi \ra 0} \p (\xi)=\infty$, $\lim _{\xi \ra \infty} \p (\xi)=0$. (See Figure $1$.)
\end{thm}

\pf
To prove the existence of solutions, we embed (\ref{l1}) into a family of problems
\beq
\lbl{l2}
u''(t)+c u'(t)+k\frac{1}{u^p(t)}=\mu +e(t), \s\s u(t+T)=u(t) \,,
\eeq
with the parameter $0 \leq k \leq 1$. When $k=0$ and $\mu=0$, we have a solution with any average $\xi$, by Lemma \ref{lma:1}. We now continue solutions with fixed average $\xi$, using Theorem \ref{thm:1}, obtaining a curve $(\xi+U(t),\mu)(k)$. By Lemma \ref{lma:6}, we have a $H^2$ bound on $U$, and then  a $H^2$ bound on $u=\xi+U$. If we choose $\xi$ large enough, we have $u(t)=\xi+U(t)>a_0>0$ for all $k$ and $t$, and hence $\frac{1}{u^p(t)}<M$, and by Lemma \ref{lma:i6} we have  a bound on $|\mu|$, so that we may continue the curve $(u(t),\mu)(k)$ for all $0 \leq k \leq 1$.
We thus have a solution of (\ref{l1}), $\hat u=\hat \xi+\hat U$  with a sufficiently large $\hat \xi$, at some $\mu =\hat \mu$.
\medskip

We now continue this solution in $\xi$. By Theorem \ref{thm:11}, we can continue locally in $\xi$. We claim that $\mu '(\xi ) \ne  0$, for all $\xi$.
Indeed, we differentiate the equation (\ref{l1}) in $\xi$
\beq
\lbl{l3}
u_{\xi}''+c u_{\xi}'+g'(u)u_{\xi}=\mu '(\xi) \,,
\eeq
with $g(u)=\frac{1}{u^p}$.
Assume  that we have $\mu '(\xi _0)=0$, at some $\xi _0$. We  set $\xi =\xi _0$ in (\ref{l3}). Then $u_{\xi}$ satisfies the equation (\ref{2}), with $h(t)=g'(u)<0$. By Lemma \ref{lma:3}, $u_{\xi}(t) \equiv 0 $. But $u(t)=\xi+U(t)$, $u_{\xi}(t)=1+U_{\xi}(t)$, so that $U_{\xi}(t)=-1$, while $U_{\xi}(t)$ is a function of average zero, a contradiction. Integrating (\ref{l1}), we have
\beq
\lbl{l4}
\mu (\xi)=\int_0^T \frac{1}{\left(\xi+U(t) \right)^p} \,dt \,.
\eeq
Increasing $\hat \xi$, if necessary, we have a uniform lower bound on $u=\xi+U$ (for $\xi>\hat \xi$), by Lemma \ref{lma:6}, and hence  the solution curve continues for all $\xi >\hat \xi$. We then see from (\ref{l4}) that $\mu (\xi) \ra 0$ as $\xi \ra \infty$. Since $\mu '(\xi ) \ne  0$, we conclude that $\mu '(\xi ) <  0$ for all $\xi$.
\medskip

We now continue the solution $(\hat \xi, \hat \mu)$ for decreasing $\xi$. Let $\xi _0 \geq 0$ be the infimum of $\xi$'s for which this continuation is possible. We claim that $\mu (\xi) \ra \infty$ as $\xi \ra \xi _0$. Indeed, assuming that $\mu$ is bounded, we see that solutions are bounded below by a positive constant, according to Lemma \ref{lma:l0}. Then the term $\frac{1}{u^p}$ is  bounded in the equation (\ref{l1}). Passing to the limit along some sequence $\xi _i \ra \xi _0$,
we obtain a bounded from below solution at $(\xi _0,\mu _0)$, which can be continued for decreasing $\xi$, contradicting  the definition of $\xi _0$.
\medskip

Finally, we show that $\xi _0=0$. Let $m=m(\xi)$, and $M=M(\xi)$ be the minimum and maximum values of $u(t)$, achieved at some points $t_0$ and $t_1$ respectively, $t_0<t_1$, ($m=u(t_0)$ and  $M=u(t_1)$). Assume, on the contrary, that $\xi _0>0$ (i.e., $\lim _{\xi \ra \xi _0} \mu(\xi)=\infty$). Then $M \geq \xi >\xi _0$ for all $\xi$, while we conclude from (\ref{l4}) that $\lim _{\xi \ra \xi _0} m(\xi)=0$. We multiply the equation (\ref{l1}) by $u'$, and integrate over $(t_0,t_1)$ (for $p \ne 1$)
\beq
\lbl{l5}
\s\s c\int _{t_0}^{t_1} {u'}^2 \, dt+\frac{1}{1-p} \left(M^{1-p}-m^{1-p} \right)=\mu (M-m)+\int _{t_0}^{t_1} e(t)u' \, dt \,.
\eeq
Denoting by $R$ the right hand side of (\ref{l5}), we estimate, using Lemma \ref{lma:6},
\beq
\lbl{l6}
R >\frac12 M \mu-c_1 \,, \s \s \mbox{for $\xi$ near $\xi _0$} \,,
\eeq
with some positive constant $c_1$. 
\medskip

\noindent
{\bf Case 1.} $\s 0<p<1$.
If $L$ is the left hand side of (\ref{l5}), we estimate ($c_2>0$)
\[
L \leq c_2+ \frac{1}{1-p} M^{1-p} \,.
\]
Combining this with (\ref{l6})
\beq
\lbl{16.1}
\frac12 M \mu-c_1< c_2+ \frac{1}{1-p} M^{1-p} \,,
\eeq
which results in a contradiction as $\mu \ra \infty$ ($M$ is bounded, since we have a uniform estimate of $|u'(t)|$, by lemma \ref{lma:6}).
\medskip

\noindent
{\bf Case 2.} $\s p=1$. Making the corresponding adjustment in (\ref{l5}), in place of (\ref{16.1}) we obtain
\[
\frac12 M \mu-c_1< c_2+  \ln M \,,
\]
resulting in a contradiction as $\mu \ra \infty$.
\medskip

\noindent
{\bf Case 3.} $\s p>1$. Then (with a constant $c_3>0$)
\beq
\lbl{l7}
L \leq c_3+ \frac{1}{p-1} m^{1-p} \,.
\eeq
Evaluating the equation (\ref{l1}) at $t_0$, we conclude that $\mu>m^{-p}-c_4$ (using that  $u''(t_0) \geq 0$), and then from (\ref{l6})
\[
R>c_4m^{-p}-c_5 \,,
\]
with positive constants $c_4$ and $c_5$.
Combining this with (\ref{l7}), and multiplying by $m^p$
\[
c_4-c_5m^{p}<c_3m^{p} + \frac{1}{p-1} m \,,
\]
which is a contradiction, since $\lim _{\xi \ra \xi _0} m(\xi)=0$.
\epf

We now turn to the case when $e(t) \in L^2$. We consider strong solutions $u(t) \in H_T^2$ of (\ref{l1}). The following two lemmas are due to R. Hakl and M. Zamora \cite{H}, who considered the $c=0$ case. Our presentation is a little different.
\begin{lma}\lbl{lma:l1} (From \cite{H})
Assume that $h(t) \in L^2(0,T)$, $u(t) \in H^2(0,T)$, both functions are defined for all $t$, and are periodic of period $T$, $u(t)>0$ for all $t$, and $c \in R$. Assume that
\beq
\lbl{l8}
u''+cu' \leq h(t) \,, \s \mbox{ for all $t$} \,.
\eeq
Then there is a positive constant $c_0$, independent of $u(t)$ and $h(t)$, such that
\beq
\lbl{l8a}
|u'(t)| \leq c_0 ||h||_{L^2(0,T)}^{2/3} \, u^{1/3}(t)   \,, \s \mbox{ for all $t$}\,.
\eeq
\end{lma}

\pf
From (\ref{l8})
\beq
\lbl{l9}
\left(e^{ct}u'(t) \right)' \leq e^{ct} |h(t)| \,.
\eeq
Let $t_0 \in [0,T]$ be an arbitrary point with $u'(t_0) \ne 0$.
\medskip

\noindent
{\em Case 1.} $u'(t_0) > 0$. Let $s_0 \in [-T,T]$ be the closest point  to the left of $t_0$ such that $u'(s_0)=0$ and $u'(t)>0$ on $(s_0,t_0)$. Multiply (\ref{l9}) by 
$\left(e^{ct}u'(t) \right)^{1/2}$, and integrate over $(s_0,t_0)$
\[
\frac23 e^{\frac32 c t_0} {u'}^{\frac32}(t_0) \leq \int_{s_0}^{t_0} e^{\frac32 c t} {u'}^{\frac12} |h(t)| \, dt \leq c_1 \left(\int_{s_0}^{t_0} u' \, dt \right)^{\frac12} ||h||_{L^2} < c_1 {u}^{\frac12}(t_0) ||h||_{L^2} \,,
\]
($c_1=e^{\frac32 c t_0}$) and (\ref{l8a}) follows.
\medskip

\noindent
{\em Case 2.} $u'(t_0) < 0$. Let $s_0 \in [0,2T]$ be the closest point  to the right of $t_0$ such that $u'(s_0)=0$ and $u'(t)<0$ on $(t_0,s_0)$.
 Multiply (\ref{l9}) by 
$\left(e^{ct}(-u') \right)^{1/2}$, and integrate over $(t_0,s_0)$
\[
\frac23 e^{\frac32 c t_0} {|u'(t_0)|}^{\frac32} \leq \int_{t_0}^{s_0} e^{\frac32 c t} {|u'|}^{\frac12} |h(t)| \, dt \,,
\]
and (\ref{l8a}) follows the same way as in Case 1.
\epf

\begin{lma}\lbl{lma:l2} (From \cite{H})
Assume that $e(t) \in \bar L_T^2$, $p \geq \frac13$, and $ \mu >0$. There is an $\ep =\ep(\mu) >0$, so that any solution $u(t) \in H^2_T$ of (\ref{l1}) satisfies $u(t)>\ep$, for all $t$.
\end{lma}

\pf
We denote $g(u)=\frac{1}{u^p}$, and $h(t)=\mu +e(t) \in L^2[0,T]$. Multiplying both sides of (\ref{l1}) by $g(u)$, and integrating, we conclude that
\beq
\lbl{l10}
\int _0^T g^2(u(t)) \, dt \leq \int _0^T h^2(t) \, dt \,.
\eeq

Integrating (\ref{l1}) over $(0,T)$, we conclude that
 $M>\delta>0$, for some $\delta=\delta(\mu)$. Let  $t_0$ and $t_1>t_0$ be the points of absolute minimum and absolute  maximum  of $u(t)$ respectively, with $u(t_0)=m$, and $u(t_1)=M$.  We now multiply both sides of (\ref{l1}) by $g^2(u)u'$, and integrate over $(t_0,t_1)$. Since the first term
\[
\int_{t_0}^{t_1} g^2u'u'' \,dt=\int_{t_0}^{t_1} g^2\left(\frac12 {u'}^2 \right)' \,dt=-\int_{t_0}^{t_1} gg'{u'}^2 \,dt \geq 0 \,,
\]
we conclude that 
\beq
\lbl{l11}
\int_{t_0}^{t_1} \frac{1}{u^{3p}} u' \,dt \leq |c| \int_{t_0}^{t_1} g^2  {u'}^2 \,dt+\int_{t_0}^{t_1} h(t)g^2u'\,dt \,.
\eeq	
We assume for the rest of the proof that $p>\frac13$, and the case $p=\frac13$ is similar. We have 
\[
\int_{t_0}^{t_1} \frac{1}{u^{3p}} u' \,dt=\frac{m^{-3p+1}}{3p-1}-\frac{M^{-3p+1}}{3p-1} \geq \frac{m^{-3p+1}}{3p-1}-\frac{{\delta}^{-3p+1}}{3p-1} \,.
\]
We now estimate the terms on the right in (\ref{l11}). By Lemma \ref{lma:l1} we have 
\beq
\lbl{l12}
\s\s \int_{t_0}^{t_1} g^2  {u'}^2 \,dt \leq c_1\int_{t_0}^{t_1} \frac{u^{2/3}}{u^{2p}} \,dt=c_1\int_{t_0}^{t_1}u^{\frac23 (-3p+1)}\,dt \leq c_2 m^{\frac23 (-3p+1)} \,,
\eeq
with some positive constants $c_1$ and $c_2$. Using (\ref{l10}),  (\ref{l12}) and Lemma \ref{lma:l1}, we estimate the last term in (\ref{l11})
\[
\int_{t_0}^{t_1} h(t)g^2u'\,dt \leq  \left( \int_{t_0}^{t_1} g^2 \,dt \right)^{1/2} \left(\int_{t_0}^{t_1} h^2 g^2  {u'}^2 \,dt\right)^{1/2} 
\]
\[
\leq c_3 \left( \int_{t_0}^{t_1}u^{\frac23 (-3p+1)} h^2\,dt \right) ^{1/2} \leq c_3 m^{\frac13 (-3p+1)} \left( \int_{t_0}^{t_1} h^2 \, dt \right)^{1/2} \leq c_4 m^{\frac13 (-3p+1)} \,,
\]
with some positive constants $c_3$ and $c_4$. Then from (\ref{l11})
\[
\frac{m^{-3p+1}}{3p-1}-\frac{{\delta}^{-3p+1}}{3p-1}\leq c_2  m^{\frac23 (-3p+1)} +c_4 m^{\frac13 (-3p+1)}\,,
\]
which results in a contradiction if $m$ is small.
\epf

We apply these lemmas to the case of $e(t) \in L^2$. 
\begin{thm}\lbl{thm:l2}
Assume that the function $e(t) \in L^2(0,T)$ is T-periodic, with $\int_0^T e(t) \, dt=0$, $p \geq \frac13$,  and $c \in R$.
Then the problem (\ref{l1}) has a unique positive $T$-periodic solution $u(t) \in W^{2,2}(0,T)$ if and only if $\mu >0$. Moreover, if we decompose $u(t)=\xi+U(t)$ with $\int_0^T U(t) \, dt=0$,  then the value of $\xi$ is a global parameter, uniquely identifying the solution pair $(u(t),\mu)$, and all positive $T$-periodic solutions lie on a unique hyperbola-like curve $\mu=\p (\xi)$. The function $\p (\xi)$ is continuous and decreasing, it is defined for  $\xi \in (\xi _0,\infty)$, with $\xi _0 \geq 0$,  and $\lim _{\xi \ra \xi _0} \p (\xi)=\infty$, $\lim _{\xi \ra \infty} \p (\xi)=0$. If, moreover, $p \leq 1$ (i.e., $\frac13 \leq p \leq 1$),  then $\xi _0=0$.
\end{thm}

\pf
The proof is similar to that of Theorem \ref{thm:l1}. We use Lemma \ref{lma:l2} instead of Lemma \ref{lma:l0} to show that solutions are bounded from below for bounded $\mu$. When it comes to the proof that $\xi _0=0$, the case $p>1$ does not carry over, since it required classical solutions (while the case $p \leq 1$ is as before).
\epf

\noindent
{\bf Remark} This result is not true if $p<1/3$. In fact, R. Hakl and M. Zamora \cite{H} show that in that case  there are $\mu$ and $e(t) \in L^2$, for which the problem (\ref{l1}) has no solution. They also considered the case $e(t) \in L^q$, where the critical $q$ turned out to be $q=\frac{1}{2p-1}$.
\medskip

\noindent
{\bf Remark} A.C. Lazer and S. Solimini \cite{LS} also considered $T$-periodic solutions of (in case $c=0$)
\[
u''+cu'-\frac{1}{u^p}=\mu+e(t) \,.
\]
When $\xi$ is small, the condition $g'(u) <\omega ^2$ is violated, and our results do not apply, although numerical computations and the result of \cite{LS} indicate that the picture is similar. When $\xi$ is large, we can proceed as before. In case of indefinite weight, the situation is more complicated, see A.J. Ure\~{n}a \cite{U}.
\medskip

\noindent
{\bf Example}  We have solved the problem (\ref{p1}) with $T =1.2$, $c=0.5$, $p=\frac 12$, $e(t)=6 \sin \frac{2 \pi}{T}  t$.
The curve $\mu=\mu (\xi)$ is given in Figure $1$. 

\begin{figure}
\begin{center}
\scalebox{0.95}{\includegraphics{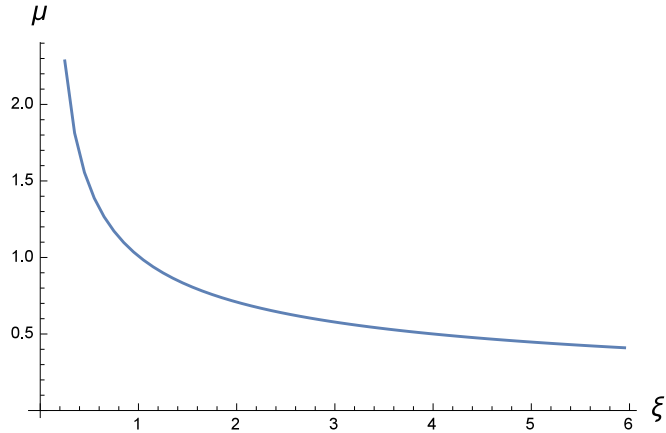}}
\end{center}
\caption{ An example for Theorem \ref{thm:l1}}
\end{figure}

\section{Electrostatic MEMS equation}
\setcounter{equation}{0}
\setcounter{thm}{0}
\setcounter{lma}{0}

Following  A. Guti\'{e}rrez and P.J. Torres \cite{G}, we now  consider an idealized mass-spring model for MEMS:
\beq
\lbl{m1}
u''(t)+c u'(t)+bu(t)+\frac{a(t)}{u^p(t)}=\mu +e(t), \s\s u(t+T)=u(t) \,,
\eeq
with the case of $p=2$, and $e(t)=0$ being physically significant. 
\begin{thm}\lbl{thm:m1}
Assume that the functions $a(t), e(t) \in C(R)$ are T-periodic, with $a(t)>0$ for all $t$, and $\int_0^T e(t) \, dt=0$; $p>0$, $c \in R$, and ($\omega = \frac{2 \pi}{T}$)
\[
0<b<\omega^2 \,.
\]
Then there is a $\mu _0>0$, so that the problem (\ref{m1}) has exactly two, one, or zero    positive $T$-periodic solutions  depending on whether  $\mu >\mu _0$, $\mu =\mu _0$, or $\mu <\mu _0$ respectively. Moreover, if we decompose $u(t)=\xi+U(t)$ with $\int_0^T U(t) \, dt=0$, then the value of $\xi$ is a global parameter, uniquely identifying the solution pair $(u(t),\mu)$, and all positive $T$-periodic solutions of (\ref{m1}) lie on a unique parabola-like curve $\mu=\p (\xi)$, defined for $\xi \in (\bar \xi,\infty)$, with some $\bar \xi >0$,  and $\lim _{\xi \ra \bar \xi} \p (\xi)=\infty$, $\lim _{\xi \ra \infty} \p (\xi)=\infty$. In case $a(t)$ is a constant, $\bar \xi =0$. (See Figure $2$.) 
\end{thm}

\pf
To prove the existence of solutions, we embed (\ref{m1}) into a family of problems
\beq
\lbl{m2}
\s\s u''(t)+c u'(t)+k\left(bu(t)+\frac{a(t)}{u^p(t)} \right)=\mu +e(t), \s\s u(t+T)=u(t) \,,
\eeq
with the parameter $0 \leq k \leq 1$. When $k=0$ and $\mu=0$, we have a solution of any average $\xi$, by Lemma \ref{lma:1}. We now continue solutions with fixed average $\xi$, using Theorem \ref{thm:1}, obtaining a curve $(\xi+U(t),\mu)(k)$. By Lemma \ref{lma:6} and Sobolev's embedding we have a uniform bound on $U$. If we choose $\xi$ large enough, we have $u(t)=\xi+U(t)>a_0>0$ for all $k$ and $t$, and hence $bu(t)+\frac{1}{u^p(t)}<M$, giving us a bound on $\mu$ by Lemma \ref{lma:i6}, so that we may continue the curve $(u(t),\mu)(k)$ for all $0 \leq k \leq 1$.
We thus have a solution of (\ref{m1}), $\hat u=\hat \xi+\hat U$  with sufficiently large $\hat \xi$, at some $\mu =\hat \mu$.
\medskip

We now continue this solution in $\xi$. By Theorem \ref{thm:11} we can continue locally in $\xi$. Increasing $\hat \xi$, if necessary, we have a uniform lower bound on $u=\xi+U$ (for $\xi>\hat \xi$), by Lemma \ref{lma:6}, and hence  the solution curve continues for all $\xi >\hat \xi$. 
Let $\bar \xi \geq 0$ be the infimum of $\xi$'s for which this continuation is possible for decreasing $\xi$. We claim that $\mu (\xi) \ra \infty$ as $\xi \ra \bar \xi$. Indeed, as in Lemma \ref{lma:l0} we show that $u(t) \geq \ep >0$, if $\mu \leq \mu _1$. By Lemma \ref{lma:6}, $U(t)$ and hence $u(t)=\xi+U(t)$ is bounded. It follows that $bu(t)+\frac{a(t)}{u^p(t)}$ is bounded, and we show as in Theorem \ref{thm:l1} that we can continue the solution curve to the left of $\bar \xi$, in contradiction with the definition of  $\bar \xi$.
\medskip

Finally, we show that $\mu (\xi)$ has only one critical point on $(\bar \xi,\infty)$, the point of global minimum. We differentiate the equation (\ref{m1}) in $\xi$
\beq
\lbl{m3}
u_{\xi}''+c u_{\xi}'+h(t)u_{\xi}=\mu '(\xi) \,,
\eeq
with $h(t)=b-p \frac{a(t)}{u^{p+1}(t)}<\omega ^2$.
At a critical point $\xi  _0$ we have $\mu '(\xi _0)=0$. We now set $\xi =\xi _0$ in (\ref{m3}). Then $w(t) \equiv u_{\xi} \mid_{\xi =\xi _0}$ satisfies the linearized problem (\ref{5}), and hence by Lemma {\ref{lma:4}, $w(t)>0$.
By Lemma \ref{lma:5}, the adjoint linear problem (\ref{9}) has a non-trivial solution $v (t)>0$. 
We differentiate (\ref{m3}) in $\xi$ again, and set $\xi =\xi _0$, obtaining
\[
u_{\xi \xi}''+c u_{\xi \xi}'+h(t)u_{\xi \xi}+h_1(t)w^2=\mu ''(\xi _0) \,,
\]
with $h_1(t)=p(p+1)\frac{a(t)}{u^{p+2}(t)}>0$.
We multiply this equation by $v(t)$ and subtract the equation (\ref{9}) multiplied by $u_{\xi \xi}$, then integrate
\[
\int _0^T h_1(t) w^2(t) v(t) \, dt=\mu ''(\xi _0) \int_0^T v(t) \, dt \,,
\]
which implies that $\mu''(\xi _0)>0$.
\epf

\medskip

\noindent
{\bf Example}  We have solved the problem (\ref{m1}) with $T=0.8$, $c=0.5$, $b=2$, $p=3$, $e(t)=5 \sin (\frac{2 \pi}{T}  t)$, $a(t)=2+ \cos ^3 (\frac{2 \pi}{T}  t)$.
The solution curve $\mu=\mu (\xi)$ is given in Figure $2$.

\begin{figure}
\begin{center}
\scalebox{0.95}{\includegraphics{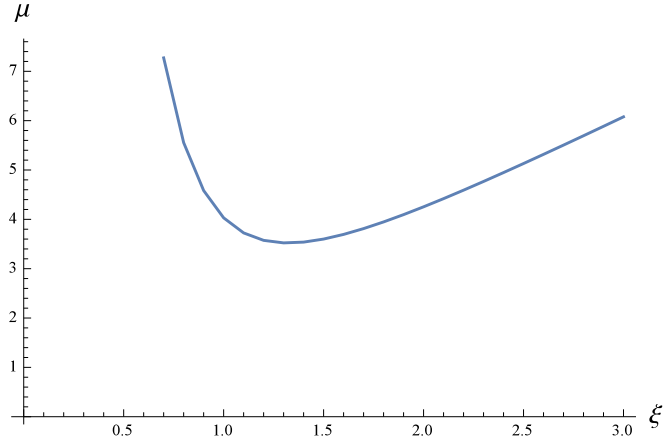}}
\end{center}
\caption{ An example for Theorem \ref{thm:m1}}
\end{figure}

\section{An equation from condensed matter physics}
\setcounter{equation}{0}
\setcounter{thm}{0}
\setcounter{lma}{0}

In the previous sections we were aided by the fact that the nonlinearities were convex.
An interesting problem in which $g(u)$ changes convexity arises as a model for fluid adsorption and wetting on a periodic corrugated substrate:
\beq
\lbl{p1}
\s u''(t)+c u'(t)+a\left(\frac{1}{u^4(t)}-\frac{1}{u^3(t)} \right)=\mu +e(t), \s\s u(t+T)=u(t) \,,
\eeq
in case $c=0$,  and $\mu =0$, see C. Rasc\'{o}n et al \cite{R}, and also a recent book of P.J. Torres \cite{T}. Similarly to G. Tarantello \cite{Ta}, we show that for suitably restricted $e(t)$, the range of $u(t)$ is not too large, so that the function $\frac{1}{u^4(t)}-\frac{1}{u^3(t)}$ is convex near the turning point, which allows a detailed description of the solution curve. We shall need the following lemma.
\begin{lma}\lbl{lma:p1}
Decompose the $T$-periodic solution of (\ref{p1}) as $u(t)=\xi+U(t)$, with $\int_0^T U(t) \, dt=0$.
Then
\[
||U(t)||_{L^{\infty}(R)} \leq \frac{\sqrt{T}}{2 \sqrt 3 |c|} ||e(t)||_{L^2[0,T]} \,.
\]
\end{lma}

\pf
We have 
\[
U''+cU'+ag(\xi +U)=\mu+e(t) \,,
\]
with $g(u)=\frac{1}{u^4}-\frac{1}{u^3}$. Multiplying by $U'$ and integrating, we get
\[
c \int_0^T {U'(t)}^2 \, dt=\int_0^T U'(t) e(t) \, dt \,,
\]
which implies that
\[
||U'(t)||_{L^2[0,T]} \leq \frac{||e(t)||_{L^2[0,T]}}{|c|} \,.
\]
The proof is completed by using the following well-known inequality for functions of zero average
\[
||U(t)||_{L^{\infty}(R)} \leq \frac{\sqrt{T}}{2 \sqrt 3 } ||U'(t)||_{L^2[0,T]} \,,
\]
see e.g., \cite{K2} for its proof.
\epf
 
\begin{thm}\lbl{thm:p1}
Assume that the function $e(t) \in C(R)$ is T-periodic, with $\int_0^T e(t) \, dt=0$,  and $a, c \in R$. Assume also that
\beq
\lbl{p2}
a \frac{3^5}{5^5}< \omega ^2 \,,
\eeq
\beq
\lbl{p3}
\frac{\sqrt{3T} ||e(t)||_{L^2[0,T]}}{|c|}< 1 \,.
\eeq
If we decompose $u(t)=\xi+U(t)$ with $\int_0^T U(t) \, dt=0$, then the value of $\xi$ is a global parameter, uniquely identifying the solution pair $(u(t),\mu)$, and all positive $T$-periodic solutions lie on a unique  continuous curve $\mu=\p (\xi)$, defined for all $\xi \in (0,\infty)$. There is a point $\xi _0 >0$, with $\mu _0=\p (\xi _0)<0$, so that the function $\p (\xi)$ is decreasing on $(0,\xi _0)$ and increasing on $(\xi _0,\infty)$, and we have $\lim _{\xi \ra 0} \p (\xi)=\infty$, $\lim _{\xi \ra \infty} \p (\xi)=0$. (See Figure $3$.)
\end{thm}

\pf
The three functions $g(u)=a \left(\frac{1}{u^4}-\frac{1}{u^3} \right)$, and its derivatives $g'(u)$ and $g''(u)$ change sign exactly once on $(0,\infty)$ at the points $u=1$, $u=4/3 $, and $u=5/3$ respectively. Then $\max _{(0,\infty)} g'(u)=g'(5/3)=a \frac{3^5}{5^5}$. By our condition (\ref{p2}), $g'(u)<\omega^2$. As before, we obtain an initial point $(\hat \xi, \hat \mu)$ on the solution curve, with $\hat \xi$ sufficiently large, and  by Theorem \ref{thm:11} we can continue solutions locally in $\xi$, and solutions continue for all $\xi >\hat \xi$, since we have a uniform bound from below. As in the proof of Theorem \ref{thm:l1}, we show that the solution curve continues to the left of $\hat \xi$ for all $\xi >0$.
\medskip

We show next that $\mu =\p (\xi)$ has only one critical point on $(0,\infty)$, the point of global minimum. Let $\xi _0$ be a critical point, $\mu '(\xi _0)=0$. We differentiate the equation (\ref{p1}) in $\xi$, and   set $\xi =\xi _0$. Then $w(t) \equiv u_{\xi} \mid_{\xi =\xi _0}$ satisfies 
\beq
\lbl{p6}
w''+cw'+g'(u)w=0 \,,
\eeq
and hence by Lemma {\ref{lma:4}, $w(t)>0$.
By Lemma \ref{lma:5}, the adjoint linear problem 
\beq
\lbl{p7a}
v''-cv'+g'(u)v=0 \,,
\eeq
has a non-trivial solution $v (t)>0$. 
We differentiate (\ref{p1}) in $\xi$ again, and set $\xi =\xi _0$, obtaining
\[
u_{\xi \xi}''+c u_{\xi \xi}'+g'(u)u_{\xi \xi}+g''(u)w^2=\mu ''(\xi _0) \,.
\]
We multiply this equation by $v(t)$ and subtract the equation (\ref{p7a}) multiplied by $u_{\xi \xi}$, then integrate
\[
\int _0^T g''(u(t)) w^2(t) v(t) \, dt=\mu ''(\xi _0) \int_0^T v(t) \, dt \,.
\]
We claim that at $\xi _0$, we have $g''(u(t))>0$, which will imply that $\mu''(\xi _0)>0$, proving that $\xi _0$ is the point of global minimum. Indeed, writing $u(t)=\xi _0+U(t)$, we have by Lemma \ref{lma:p1}, and the condition (\ref{p3}),
\beq
\lbl{p7}
|U(t)|<\frac16 \,, \s \mbox{for all $t$} \,.
\eeq
Integrating (\ref{p6}), we get
\[
\int _0^T g'(u(t))w \, dt=0 \,.
\]
Hence, $u(t)$ takes on values where $g'(u)$ is both positive and negative, i.e., the range of $u(t)$ includes $u=4/3$. By (\ref{p7}), $\xi _0<4/3+1/6$, and then $u(t)=\xi _0+U(t)<4/3+1/6+1/6=5/3$, and hence $g''(u(t))>0$.
\medskip

Since $\mu _0=\mu (\xi _0)$ is  the absolute minimum of $\mu(\xi)$, to show that $\mu (\xi _0)<0$ it suffices to show that $\mu(\xi)$ takes on negative values. Integrating (\ref{p1}), we get
\beq
\lbl{p8}
\mu=\frac{a}{T} \int_0^T \frac{1-\xi-U(t)}{(\xi+U(t))^4} \, dt<0 \,,
\eeq
for $\xi$ large, in view of the uniform estimate of $U$ given by (\ref{p7}).
\epf

In the physically significant case when $c=0$,  and $\mu =0$, this theorem implies the existence  and uniqueness of a $T$-periodic solution. Moreover, the average value of this solution is greater than $1$, i.e., $\xi >1$. Indeed, from (\ref{p8})
\[
\mu(1)=-\frac{a}{T} \int_0^T \frac{U(t)}{(1+U(t))^4} \, dt>0 \,,
\]
(recall that $\int_0^T U(t) \, dt=0$) so that $\xi >1$ when $\mu =0$.
\medskip

\noindent
{\bf Example}  We have solved the problem (\ref{p1}) with $c=0.3$, $a=3$, $e(t)=8 \cos2 \pi  t$ (i.e., $T =1$).
The curve $\mu=\mu (\xi)$ is given in Figure $3$. 

\begin{figure}
\begin{center}
\scalebox{0.95}{\includegraphics{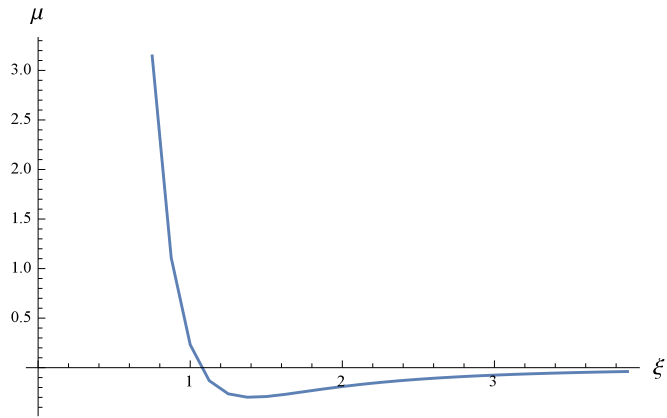}}
\end{center}
\caption{ An example for Theorem \ref{thm:p1}}
\end{figure}

\section{Numerical computation of solutions}
\setcounter{equation}{0}
\setcounter{thm}{0}
\setcounter{lma}{0}

To  solve the periodic problem (with $g(t+T,u)=g(t,u)$, for all $t$ and $u$)
\beq
\lbl{n1}
u''+cu'+g(t,u)=\mu+e(t) \,, \; u(t)=u(t+T), \;\; u'(t)=u'(t+T)
\eeq
we used continuation in $\xi$. We began by implementing the numerical solution of the following linear periodic problem: given the $T$-periodic functions  $b(t)$ and $f(t)$, and a constant $c$, find the $T$-periodic solution of
\beq
\lbl{n2}
\s\s\s L[y] \equiv y''(t)+c y' +b(t)y=f(t), \;\; y(t)=y(t+T), \; y'(t)=y'(t+T) \,.
\eeq
The general solution of (\ref{n2}) is of course
\[
y(t)=Y(t)+c_1 y_1(t)+c_2 y_2(t) \,,
\]
where $Y(t)$ is a particular solution, and $y_1$, $y_2$ are two solutions of the corresponding homogeneous equation. To find $Y(t)$, we used the NDSolve command to solve (\ref{n2}) with $y(0)=0$, $y'(0)=1$. {\em Mathematica} not only solves differential equations numerically, but it returns the solution as an interpolated function of $t$, practically indistinguishable from an explicitly defined solution. (We believe that the NDSolve command is a {\em game-changer}, making equations with variable coefficients as easy to handle, as the ones with constant coefficients.)  We calculated $y_1$ and  $y_2$ by solving the corresponding homogeneous equation with the initial conditions $y_1(0)=0$, $y_1'(0)=1$, and $y_2(0)=1$, $y_2'(0)=0$. We then select $c_1$ and $c_2$, so that
\[
y(0)=y(T), \s y'(0)=y'(T) \,,
\]
which is just a linear $2 \times 2$ system. This gives us the $T$-periodic solution of (\ref{n2}), or $L^{-1}[f(t)]$, where $L[y]$ denotes the left hand side of (\ref{n2}), subject to the periodic boundary conditions. 
\medskip

Then we have implemented the {\em ``linear solver"}, i.e., the numerical solution of the following problem: given the $T$-periodic functions  $b(t)$ and $f(t)$, and a constant $c$, find the constant $\mu$, so that the problem
\beq
\lbl{n10}
 y''(t) +c y' +b(t)y=\mu+f(t), \; \int_0^T y(t) \, dt=0 
\eeq
has a $T$-periodic solution of zero average, and compute that solution $y(t)$. The solution is 
\[
y(t)=L^{-1}[f(t)]+\mu L^{-1}[1] \,,
\]
with the constant  $\mu$ chosen so that $\int_0^T y(t) \, dt=0$.
\medskip

Turning to the problem (\ref{n1}), we begin with an initial $\xi _0$, and using a step size $\Delta \xi$, we compute the solution $(\mu _n ,u_n(t))$ with the average of  $u_n(t) $ equal to  $\xi _n=\xi _0 +n \Delta \xi$, $n=1,2, \ldots, nsteps$, in the form $u_n(t)=\xi _n+U_n(t)$, where $U_n(t)$ is the $T$-periodic solution of
\beq
\lbl{n3}
\s\s\s U''+cU'+g(t,\xi _n+U)=\mu+e(t) \,, \; \int_0^T U(t) \, dt=0 \,.
\eeq
With $U_{n-1}(t)$  already computed, we use Newton's method to solve for $U_{n}(t)$. We linearize the equation (\ref{n3}) at $\xi_{n}+U_{n-1}(t)$, writing $g(t,\xi _n+U)=g(t,\xi _n+U_{n-1})+g_u(t,\xi _{n}+U_{n-1})(U-U_{n-1})$, and call  the linear solver to find the $T$-periodic solution of the problem (\ref{n10}) with $b(t)=g_u(t,\xi _n+U_{n-1})$, and $f(t)=e(t)-g(t,\xi _{n}+U_{n-1})+g_u(t,\xi _{n}+U_{n-1})U_{n-1}$, obtaining an approximation of $U_n$ and $\mu _n$, call them $\bar U_n$ and $\bar \mu _n$. We then linearize the equation (\ref{n3}) at $\xi_{n}+\bar U_{n}(t)$, to get a better approximation of $U_n$ and $\mu _n$, and so on, making several iterative steps. (In short, we use Newton's method to solve (\ref{n3}), with $U_{n-1}$ being the initial guess.)
We found that just two iterations of Newton's method, coupled with a relatively small $\Delta \xi$ (e.g., $\Delta \xi=0.1$), were sufficient for accurate computation of the solution curves. Finally, we plot the points $(\xi _n, \mu _n)$ to obtain the solution curve.
\medskip

We have verified our numerical results by an independent calculation. Once  a periodic solution $u(t)$ is computed at some $\mu$, we took its values of $u(0)$ and $u'(0)$, and computed numerically the solution of (\ref{n1}) at this $\mu$, with these initial data (using the NDSolve command), as well as the  average of $u(t)$. We had a complete agreement (with $u(t)$ and $\xi$) for all $\mu$, and all equations that we tried.

\end{document}